\documentclass[12pt]{article}
\usepackage{a4}
\usepackage{subfig,graphicx}
\usepackage{amsthm,amssymb,amsmath}
\usepackage{enumitem}
\usepackage{xcolor}
\usepackage{eucal}

\newtheorem{theorem}{Theorem}

\newtheorem{corollary}[theorem]{Corollary}
\newtheorem{conjecture}[theorem]{Conjecture}

\newcommand\size[1] {\left|{#1}\right|}
\newcommand\Set[2] {\left\{{#1}:\,{#2}\right\}}
\newcommand\Setx[1] {\left\{{#1}\right\}}

\newcommand{\FF}{\mathbb F}
\newcommand{\QQ}{\mathbb Q}
\newcommand{\RR}{\mathbb R}
\newcommand{\TT}{\mathcal R}
\newcommand{\ZZ}{\mathbb Z}
\newcommand{\fln}{\Phi}
\newcommand{\flnc}{\Phi_c}
\newcommand{\flnm}{\Phi_\mathrm{mod}}
\newcommand{\flns}{\Phi_\mathrm{z}}
\newcommand{\ispec}[1]{\overline{\mathcal{S}}(#1)}
\newcommand{\mspec}[1]{\mathcal{S}_\mathrm{mod}(#1)}
\newcommand{\spec}[1]{\mathcal{S}(#1)}
\newcommand{\sig}{\Sigma}
\newcommand{\Macajova}{M\'{a}\v{c}ajov\'{a}}
\newcommand{\Skoviera}{\v{S}koviera}

\newcommand{\fig}[1]{{%
    \includegraphics[page=#1]{signed-figs}}}
\newcommand{\hf}{\hspace*{0mm}\hfill\hspace*{0mm}}

\usepackage[utf8]{inputenc}
\usepackage[
  style=numeric-comp, 
  url=false,
  doi=false,
  eprint=true,
  sorting=nyt,
  isbn=false,
  firstinits=true,
  defernumbers=true,
  labelnumber=true,
  maxnames=10,
  backend=bibtex]{biblatex} 
\usepackage{xpatch}
\xpatchbibmacro{volume+number+eid}{%
  \setunit*{\adddot}%
}{%
}{}{}
\DeclareFieldFormat[article,thesis,inbook,inproceedings]{title}{#1} 
\DeclareFieldFormat[article,thesis,inbook,inproceedings]{pages}{#1} 
\DeclareFieldFormat[article]{number}{~\mkbibparens{#1}}
\DeclareFieldFormat[article]{volume}{\textbf{#1}}
\renewbibmacro{in:}{%
  \ifentrytype{article}{}{\printtext{\bibstring{in}\intitlepunct}}}
\AtEveryBibitem{\clearfield{month}}
\AtEveryBibitem{\clearfield{day}}

\title{\textbf{Nowhere-zero flows in signed graphs:\\A survey}}%
\author{Tom\'{a}\v{s} Kaiser\thanks{Department of Mathematics,
    Institute for Theoretical Computer Science (CE-ITI), and European
    Centre of Excellence NTIS (New Technologies for the Information
    Society), University of West Bohemia, Pilsen, Czech
    Republic. Supported by project GA14-19503S of the Czech Science
    Foundation. Email: \texttt{kaisert@kma.zcu.cz}.}  \and Robert
  Lukot'ka\thanks{Department od Computer Science, Faculty of Mathematics, 
  Physics, and Informatics, Comenius University, Bratislava, Slovakia. 
  Supported by project VEGA 1/0876/16. Email: \texttt{lukotka@dcs.fmph.uniba.sk}.} 
  \and Edita Rollov\'{a}  \thanks{European
    Centre of Excellence NTIS (New Technologies for the Information
    Society), University of West Bohemia, Pilsen, Czech
    Republic. Partially supported by project GA14-19503S of the Czech Science Foundation and by project LO1506 of the Czech Ministry of Education, Youth and Sports. 
    Email: \texttt{rollova@ntis.zcu.cz}.}} \date{}

\begin{document}
\maketitle

\begin{abstract}
  We survey known results related to nowhere-zero flows and related
  topics, such as circuit covers and the structure of circuits of
  signed graphs. We include an overview of several different
  definitions of signed graph colouring. \medskip
	
\end{abstract}


\section{Introduction}
\label{sec:intro}

Nowhere-zero flows were defined by
Tutte~\cite{tutte_contribution_1954} as a dual problem to
vertex-colouring of (unsigned) planar graphs. Both notions have been
extended to signed graphs. The definition of nowhere-zero flows on
signed graphs naturally comes from the study of embeddings of graphs
in non-orientable surfaces, where nowhere-zero flows emerge as the
dual notion to local tensions. 
There is a close relationship between nowhere-zero flows and circuit
covers of graphs as every nowhere-zero flow on a graph $G$ determines
a covering of $G$ by circuits. This relationship is maintained for
signed graphs, although a signed version of the definition of circuit
is required.

The area of nowhere-zero flows on (signed) graphs recently received a
lot of attention thanks to the breakthrough result of
Thomassen~\cite{thomassen_weak_2012} on nowhere-zero $3$-flows in
unsigned graphs. The purpose of this paper is to capture the current
state of knowledge regarding nowhere-zero flows and circuit covers in
signed graphs.  To cover the latest developments in this active area,
we decided to include some papers that are currently still in the
review process and only available from the arXiv. Clearly, until the
papers are published, the results have to be taken with caution.

An extensive source of links to material on signed graph is the
dynamic survey of Zaslavsky~\cite{zaslavsky_mathematical_2012}. For
the theory of unsigned graphs in general, see~\cite{bondy_graph_2008}
or \cite{diestel_graph_2005}. Nowhere-zero flows and circuit covers in
unsigned graphs are the subject of Zhang~\cite{zhang_integer_1997};
see also \cite{zhang_circuit_2012}.
\bigskip

Graphs in this survey may contain parallel edges and loops. A circuit
is a connected 2-regular graph.

A \emph{signed graph} is a pair $(G,\sig)$, where $\sig$ is a subset
of the edge set of $G$. The edges in $\sig$ are \emph{negative}, the
other ones are \emph{positive}. An \emph{unbalanced circuit} in
$(G,\sig)$ is an (unsigned) circuit in $G$ that has an odd number of
negative edges. A \emph{balanced circuit} in $(G,\sig)$ is an
(unsigned) circuit in $G$ that is not unbalanced. A subgraph of $G$ is
\emph{unbalanced} if it contains an unbalanced circuit; otherwise, it
is \emph{balanced}. A signed graph is \emph{all-negative}
(\emph{all-positive}) if all its edges are negative (positive,
respectively).

Signed graphs were introduced by Harary~\cite{harary_notion_1953} as a
model for social networks. Zaslavsky~\cite{zaslavsky_signed_1982-1}
then employed the following kind of equivalence on the class of signed
graphs. Given a signed graph $(G,\sig)$, \emph{switching} at a vertex
$v$ is the inversion of the sign of each edge incident with $v$. Thus,
the signature of the resulting graph is the symmetric difference of
$\sig$ with the set of edges incident with $v$. Signed graphs are said
to be \emph{equivalent} if one can be obtained from the other by a
series of switchings. It is easy to see that equivalent signed graphs
have the same sets of unbalanced circuits and the same sets of
balanced circuits.

Given a partition $\Setx{A,B}$ of the vertex set of a graph $G$, we
let $[A,B]$ denote the set of all edges with one end in $A$ and one in
$B$. Harary~\cite{harary_notion_1953} characterised balanced graphs:
\begin{theorem}
  A signed graph $(G,\sig)$ is balanced if and only if there is a set
  $X\subseteq V(G)$ such that the set of negative edges is precisely
  $[X,V(G)-X]$.
\end{theorem}
It follows from the characterisation that a signed graph is balanced
if it is equivalent to an all-positive signed graph. We say that a
graph is \emph{antibalanced} if it is equivalent to an all-negative
signed graph.

\begin{figure}
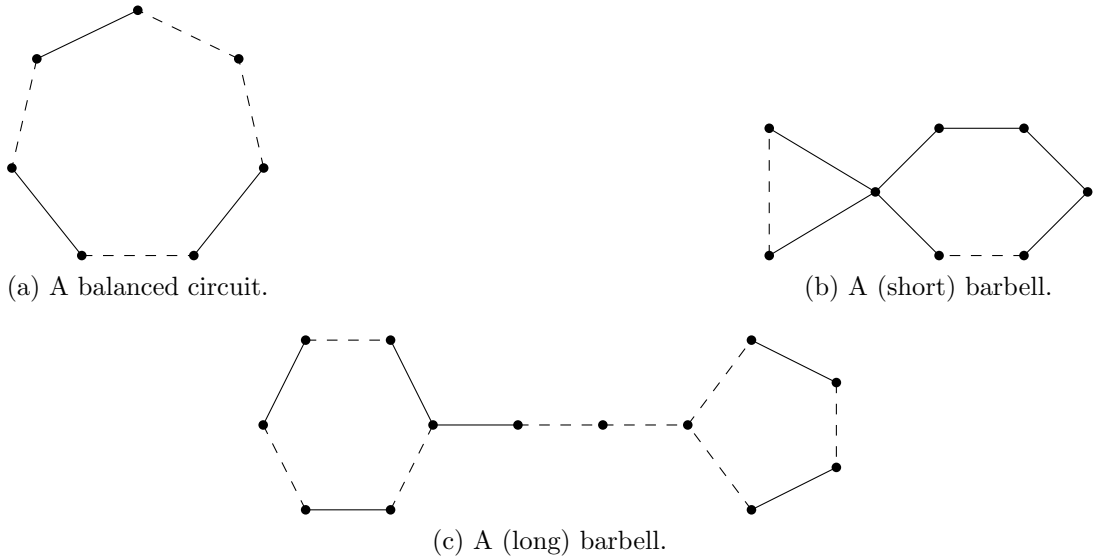

  \centering
  \subfloat[][A balanced circuit.]{\fig1}\hf
  \subfloat[][A (short) barbell.]{\fig2}\\
  \subfloat[][A (long) barbell.]{\fig3}
  \caption{Signed circuits. Dashed lines indicate negative edges.}
  \label{fig:sig-c}
\end{figure}

The role of circuits in unsigned graphs is played by signed circuits
(known to be the circuits of the associated signed graphic matroid
\cite{zaslavsky_signed_1982}). A \emph{signed circuit} belongs to one
of two types (cf. Figure~\ref{fig:sig-c}):
\begin{itemize}
\item a balanced circuit,
\item a \emph{barbell}, i.e., the union of two unbalanced circuits
  $C_1,C_2$ and a (possibly trivial) path $P$ with endvertices
  $v_1\in V(C_1)$ and $v_2\in V(C_2)$, such that $C_1-v_1$ is disjoint
  from $P\cup C_2$ and $C_2-v_2$ is disjoint from $P\cup C_1$.
\end{itemize}

If the path in a barbell is trivial, the barbell is sometimes called
\emph{short}; otherwise it is \emph{long}. If we need to speak about a
circuit in the usual unsigned sense, we call it an \emph{ordinary
  circuit}.

\begin{figure}
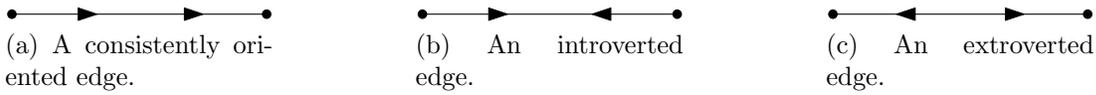

  \centering
  \subfloat[][A consistently oriented edge.]{\fig4}\hf
  \subfloat[][An introverted edge.]{\fig5}\hf
  \subfloat[][An extroverted edge.]{\fig6}
  \caption{Oriented edges of a signed graph.}
  \label{fig:edges}
\end{figure}

To obtain an orientation of a signed graph, each edge of $(G,\sig)$ is
viewed as composed of two half-edges. An \emph{orientation} of an edge
$e$ is obtained by giving each of the two half-edges $h,h'$ making up
$e$ a direction. We say that $e$ is \emph{consistently oriented} if
exactly one of $h,h'$ is directed toward its endvertex. Otherwise, it
is \emph{extroverted} if both $h$ and $h'$ point toward their
endvertices, and \emph{introverted} if none of them does. (See
Figure~\ref{fig:edges}.)  We say that an oriented edge $e$ is
\emph{incoming} at its endvertex $v$ if its half-edge incident with
$v$ is directed toward $v$, and that $e$ is \emph{outgoing} at $v$
otherwise.

An \emph{orientation} (\emph{bidirection}) of $(G,\sig)$ is the
assignment of an orientation to each edge of $G$ (in the above sense), in such a way that the positive edges are exactly the consistently oriented ones. An oriented signed graph is often called a \emph{bidirected graph}.


\section{Nowhere-zero flows}
\label{sec:nowhere-zero-flows}

\subsection{Introduction}
\label{sec:flows-intro}
Let $\Gamma$ be an Abelian group. A \emph{$\Gamma$-flow} in $(G,\sig)$
consists of an orientation of $(G,\sig)$ and a function
$\varphi:\,E(G)\to\Gamma$ such that the usual conservation law is
satisfied --- that is, for each vertex $v$, the sum of $\varphi(e)$
over the incoming edges $e$ at $v$ equals the sum of $\varphi(e)$ over
the outgoing ones. A $\ZZ$-flow is said to be a \emph{$k$-flow} (where
$k\geq 2$ is an integer) if $|\varphi(e)| < k$ for each edge $e$. A
$\Gamma$-flow is \emph{nowhere-zero} if the value $0$ is not used at
any edge. If a signed graph $(G,\sig)$ admits a nowhere-zero
$\ZZ$-flow, then its \emph{flow number} $\fln(G,\sig)$ is defined as
the least $k$ such that $(G,\sig)$ admits a nowhere-zero
$k$-flow. Otherwise, $\fln(G,\sig)$ is defined as $\infty$. 

A signed graph is said to be \emph{flow-admissible} if it admits at least one nowhere-zero $\ZZ$-flow. Bouchet~\cite{bouchet_nowhere-zero_1983} observed the following characterisation in combination with \cite{zaslavsky_signed_1982}.

\begin{theorem}\label{t:flow-admissible}
 A signed graph $(G,\sig)$ is flow-admissible if and only if every edge of $(G,\sig)$ belongs to a signed circuit of $(G,\sig)$.
\end{theorem}

Thus, an all-positive signed graph is flow-admissible if and only if
it is bridgeless. Nowhere-zero flows on signed graphs are, in fact, a generalization of the same concept on unsigned graphs, because the definition of a flow for all-positive
signed graph corresponds to the usual definition of a flow on an
unsigned graph.

This useful corollary directly follows from the previous theorem:
\begin{corollary}\label{t:one-edge}
A signed graph with one negative edge is not flow-admissible. 
\end{corollary}

For 2-edge-connected unbalanced signed graphs the converse is also true \cite{macajova_nowhere-zero_2014}.

\begin{theorem}\label{t:special_edge}
 An unbalanced 2-edge-connected signed graph is flow-admissible except in the case when it contains an edge whose removal leaves a balanced graph.
\end{theorem}

Tutte~\cite{tutte_imbedding_1949} proved that an unsigned graph $G$
admits a nowhere-zero $k$-flow if and only if it admits a nowhere-zero
$\ZZ_k$-flow. He proved, in fact, that the number of nowhere-zero
$\Gamma$-flows only depends on $\size\Gamma$ and not on the structure
of $\Gamma$. This is not true for signed graphs in general. For
instance, an unbalanced circuit admits a nowhere-zero $\ZZ_2$-flow,
but it does not admit any integer flow. It was, however, proved for
signed graphs that if $|\Gamma|$ is odd, then the number of
nowhere-zero $\Gamma$-flows does not depend on the structure of
$\Gamma$~\cite{beck_number_2006}. In the same paper, Beck and
Zaslavsky further proved that if $|\Gamma|$ is odd, the number of
nowhere-zero $\Gamma$-flows on a signed graph $(G,\sig)$ is a
polynomial in $|\Gamma|$, independent of the actual group, extending
the result of Tutte~\cite{tutte_contribution_1954} for unsigned graphs
(where arbitrary $|\Gamma|$ is allowed).

\subsection{Bouchet's conjecture}
\label{sec:bouchet}

Nowhere-zero flows in signed graphs were introduced by Edmonds and
Johnson~\cite{edmonds_matching_1970} for expressing algorithms on
matchings, but systematically studied at first by Bouchet in the
paper~\cite{bouchet_nowhere-zero_1983}. Bouchet also stated a
conjecture that parallels Tutte's $5$-Flow Conjecture and occupies a
similarly central place in the area of signed graphs:

\begin{conjecture}[Bouchet]\label{conj:bouchet}
  Every flow-admissible signed graph admits a nowhere-zero $6$-flow.
\end{conjecture}

\begin{figure}
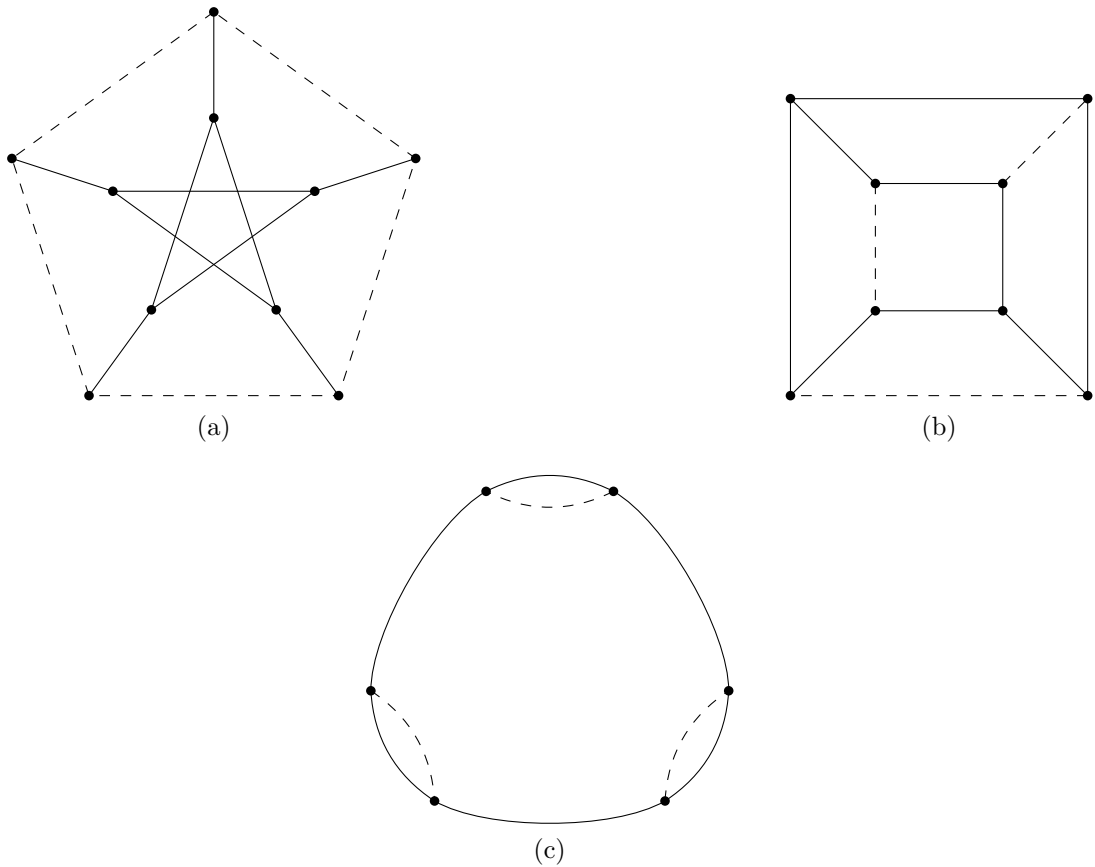

  \centering
  \subfloat[][]{\fig7}\hf
  \subfloat[][]{\fig8}\\
  \subfloat[][]{\fig9}
  \caption{Signed graphs without nowhere-zero $5$-flows. The example
    (c) gives rise to an infinite family.}
  \label{fig:bouchet}
\end{figure}

The value $6$ would be best possible since Bouchet showed that the Petersen graph with a signature corresponding to its triangular embedding in the projective plane (see Figure~\ref{fig:bouchet}(a))
admits no nowhere-zero $5$-flow. Another signed graph with this
property given in Figure~\ref{fig:bouchet}(b) was found using a computer search by M\' a\v cajov\' a~\cite{macajova_personal_2012}; an
infinite family of such graphs containing the one in
Figure~\ref{fig:bouchet}(c) was constructed by Schubert and
Steffen~\cite{schubert_nowhere-zero_2015}. (The other graphs in the family are obtained by arranging
$k$ unbalanced $2$-circuits, $k\geq 5$ odd, in a circular manner as in
Figure~\ref{fig:bouchet}(c).) Interestingly, no other examples of
signed graphs with flow number $6$ appear to be known 
(except for signed graphs constructed from the known ones by simple flow reductions).

In~\cite{bouchet_nowhere-zero_1983}, Bouchet proved that the
conjecture holds with the number $6$ replaced by $216$. The constant
was improved by Z\'{y}ka~\cite{zyka_nowhere-zero_1987} to $30$;
according to~\cite{devos_flows_2013}, the same result was obtained by
Fouquet (unpublished). Currently, the best general bound is due to DeVos~\cite{devos_flows_2013}:
\begin{theorem}[DeVos]\label{t:devos}
  Every flow-admissible signed graph admits a nowhere-zero $12$-flow.
\end{theorem}


\subsection{Higher edge-connectivity}
\label{sec:higher}

Assumptions about the connectivity of a graph usually allow for better bounds on its flow number, or at least make the bounds easier to establish. In the unsigned case, for example, it is well-known that any $4$-edge-connected graph admits a nowhere-zero $4$-flow, so in particular these graphs are not interesting with respect to Tutte's 5-Flow Conjecture.

For signed graphs, connectivity alone is not sufficient, as can be seen by Corollary~\ref{t:one-edge}: no signed graph with a single negative edge is flow-admissible. Thus, an assumption of flow-admissibility needs to be included.

Very recently, the following was shown by Cheng et al.~\cite{cheng_integer_2016} for 2-edge-connected signed graphs:

\begin{theorem}[Cheng, Lu, Luo, Zhang]\label{t:2-edge-connected}
Every flow-admissible $2$-edge-connected signed graph has a nowhere-zero $11$-flow.
\end{theorem}

The flow number of $3$-edge-connected flow-admissible signed graphs was bounded by $25$ in \cite{wei_flows_2011}, which was improved to $15$ in \cite{wei_nowhere-zero_2014}, and to $9$ in~\cite{yang_nowhere-zero_2015}. 
A further improvement, based on the ideas of the proof of the Weak $3$-Flow Conjecture~\cite{thomassen_weak_2012,lovasz_nowhere-zero_2013-1} was obtained by Wu et al.~\cite{wu_nowhere-zero_2014}:

\begin{theorem}[Wu, Ye, Zang and Zhang]\label{t:3-edge-connected}
  Every flow-admissible $3$-edge-connected signed graph admits a
  nowhere-zero $8$-flow.
\end{theorem}

Khelladi~\cite{khelladi_nowhere-zero_1987} proved already in 1987 that any flow-admissible $4$-edge-connected signed graph has a nowhere-zero $18$-flow, and established Bouchet's Conjecture for flow-admissible $3$-edge-connected graphs containing no long barbell. 
Bouchet's Conjecture was also established for $6$-edge-connected signed graphs~\cite{xu_flows_2005}.
The best possible result for $4$-edge-connected signed graphs was given by Raspaud and Zhu~\cite{raspaud_circular_2011} using signed circular flows (for further results regarding signed circular flows, see Section~\ref{sec:circular}):

\begin{theorem}[Raspaud and Zhu]\label{t:4-edge-connected}
Every flow-admissible $4$-edge-connected signed graph admits a
  nowhere-zero $4$-flow.
\end{theorem}




\subsection{Signed regular graphs}
\label{sec:regular}

Most of the research in the area of signed regular graphs is focused on signed cubic graphs. \Macajova{} and \Skoviera~\cite{macajova_remarks_2015} characterised
signed cubic graphs with flow number $3$ or $4$. 

Let $(G,\sig)$ be a signed graph and let $\Setx{X_1,X_2}$ be a partition of $V(G)$.
If set of positive edges of $G$ is exactly $[X_1,X_2]$, then the partition $\Setx{X_1,X_2}$ is \emph{antibalanced}. A signed graph $(G,\sig)$ has an antibalanced partition if and only if it is antibalanced~\cite{harary_notion_1953}. For connected antibalanced signed graphs, the antibalanced partition is uniquely determined. The circuit $(Z,\sig)$ with antibalanced partition $\Setx{X_1,X_2}$ is \emph{half-odd} if for some $i=1,2$, each component of $Z-X_i$ is either $Z$ itself or a path of odd length.

\begin{theorem}
  Let $G,\sig$ be a signed cubic graph. Then the following hold:
  \begin{enumerate}[label=(\roman*)]
  \item $(G,\sig)$ admits a nowhere-zero $3$-flow if and only if it is
    antibalanced and has a perfect matching,
  \item $(G,\sig)$ admits a nowhere-zero $4$-flow if and only if an
    equivalent signed graph contains an antibalanced $2$-factor with
    all components half-odd whose complement is an all-negative
    perfect matching.
  \end{enumerate}
\end{theorem}

Group-valued nowhere-zero flows in signed cubic graphs were also
investigated in~\cite{macajova_remarks_2015}, with the following
results:
\begin{theorem}\label{t:grp-flows}
  Let $(G,\sig)$ be a signed cubic graph. Then the following holds:
  \begin{enumerate}[label=(\roman*)]
  \item $(G,\sig)$ has a nowhere-zero $\ZZ_3$-flow if and only if it
    is antibalanced,
  \item $(G,\sig)$ has a nowhere-zero $\ZZ_4$-flow if and only if it
    has an antibalanced $2$-factor, 
  \item $(G,\sig)$ has a nowhere-zero $\ZZ_2\times\ZZ_2$-flow if and
    only if it $G$ is $3$-edge-colourable.
  \end{enumerate}
\end{theorem}

The \emph{integer flow spectrum} $\ispec G$ of a graph $G$ is the set
of all $t$ such that there is some signature $\sig$ of $G$ such that
$\fln(G,\sig) = t$. The analogous notion where $\fln$ is replaced by
$\flnc$ (circular flow number, for the definition see Section~\ref{sec:circular}) 
is called the \emph{flow spectrum} and denoted by $\spec G$.

Schubert and Steffen~\cite{schubert_nowhere-zero_2015} investigated
both kinds of flow spectra in signed cubic graphs. They proved the
following ($K^3_2$ denotes the graph on two vertices with three parallel edges and no loops):
\begin{theorem}
  If $G$ is a cubic graph different from $K^3_2$, then we have the
  following equivalences:
  \begin{equation*}
    G\text{ has a $1$-factor}\iff 3\in\spec G\iff 3\in\ispec G\iff
    4\in\ispec G.
  \end{equation*}
  Moreover, if $G$ has a $1$-factor, then $\Setx{3,4}\subseteq\ispec
  G\cap \spec G$.
\end{theorem}
Moreover, there exists a graph with $\spec G = \{3,4\}$: the cubic
graph on $4$ vertices that has two double edges; and an infinite
family of graphs with
$\spec G = \{3,4,6\}$~\cite{schubert_nowhere-zero_2015}.

For general $(2t+1)$-regular graphs the following was proved in~\cite{schubert_nowhere-zero_2015}:
\begin{theorem}
Let $(G,\sig)$ be a signed graph. If $G$ is $(2t+1)$-regular and has a 1-factor, then $3\in \ispec G$.
\end{theorem}

A rather different version of a `flow spectrum' was introduced by
\Macajova{} and \Skoviera{}~\cite{macajova_remarks_2015}. Let the
\emph{modular flow number} $\flnm(G,\sig)$ be the least integer $k$
such that $(G,\sig)$ admits a nowhere-zero $\ZZ_k$-flow. It is easy to
see that $\flnm(G,\sig) \leq \fln(G,\sig)$. The \emph{modular flow
  spectrum} $\mspec{G,\sig}$ is the set of all $k$ such that
$(G,\sig)$ admits a nowhere-zero $\ZZ_k$-flow.

It is shown in~\cite{macajova_remarks_2015} that the modular flow
spectrum of a signed cubic graph may contain gaps. Namely, a signed
cubic graph is given whose modular flow spectrum equals
$\Setx{3,5,6,7,\dots}$. On the other hand, it is also shown that
bridges are essential in this example:
\begin{theorem}
  Let $(G,\sig)$ be a bridgeless cubic signed graph with
  $\flnm(G,\sig) = 3$. Then $(G,\sig)$ admits a nowhere-zero
  $\Gamma$-flow for any Abelian group $\Gamma$ of order at least $3$,
  except possibly for $\ZZ_2\times\ZZ_2$. In particular,
  $\mspec{G,\sig} = \Set k {k\geq 3}$.
\end{theorem}


\subsection{Other classes of signed graphs}
\label{sec:special}

More detailed information about nowhere-zero flows has been obtained
for some other classes of signed graphs, namely signed Eulerian,
complete and complete bipartite graphs, series-parallel graphs, Kotzig graphs, and signed graphs with two
negative edges.

Unsigned Eulerian graphs enjoy the minimum possible flow number of
$2$. This is not the case with signed Eulerian graphs, for which the
following was proved in~\cite{xu_flows_2005}:
\begin{theorem}\label{t:xz-eulerian}
  A connected signed graph has a nowhere-zero $2$-flow if and only if
  it is Eulerian and the number of its negative edges is even.
\end{theorem}

Flows in signed Eulerian graphs were further investigated by
\Macajova{} and \Skoviera~\cite{macajova_nowhere-zero_2014} (cf. also
the extended abstract~\cite{macajova_determining_2011}). Their results
are summarised by Theorem~\ref{t:eulerian} below. Before we state it,
we introduce one more notion. A signed Eulerian graph is \emph{triply
  odd} if it has an edge-decomposition into three Eulerian subgraphs
sharing a vertex, each of which contains an odd number of negative
edges.
\begin{theorem}\label{t:eulerian}
  Let $(G,\sig)$ be a connected signed Eulerian graph. Then
  \begin{enumerate}[label=(\roman*)]
  \item if $(G,\sig)$ is flow-admissible, then it has a nowhere-zero
    $4$-flow,
  \item $\fln(G,\sig) = 3$ if and only if $(G,\sig)$ is triply odd.
  \end{enumerate}
\end{theorem}
Observe that together with Theorem~\ref{t:xz-eulerian} and Corollary~
\ref{t:one-edge}, this result implies a full characterisation of
signed Eulerian graphs of each of the possible flow numbers ($2$, $3$, $4$ and $\infty$).  

Another result of~\cite{macajova_nowhere-zero_2014} is a
characterisation, for any Abelian group $\Gamma$, of signed Eulerian
graphs admitting a nowhere-zero $\Gamma$-flow. Namely:

\begin{theorem}
Let $(G,\sig)$ be a signed eulerian graph and let $\Gamma$ be a nontrivial Abelian group.
The following statements hold true.
(a) If $\Gamma$ contains an involution, then $G$ admits a nowhere-zero $\Gamma$-flow.
(b) If $\Gamma\cong \ZZ_3$, then $(G,\sig)$ admits a nowhere-zero $\Gamma$-flow if and only if $G$ is triply odd.
(c) Otherwise, $G$ has a nowhere-zero $\Gamma$-flow if and only if $G$ is flow-admissible.
\end{theorem}
\smallskip

\Macajova{} and Rollov\'{a}~\cite{macajova_nowhere-zero_2015} (cf. also the extended abstract \cite{macajova_flow_2011}) determined the
complete characterisation of flow numbers of signed complete and signed complete bipartite graphs:
\begin{theorem}\label{t:complete}
  Let $(K_n,\sig)$ be a flow-admissible signed complete graph. Then
  \begin{enumerate}[label=(\roman*)]
  \item $\fln(K_n,\sig) = 2$ if and only if $n$ is odd and $\size\sig$
    is even,
  \item $\fln(K_n,\sig) = 4$ if and only if $(K_n,\sig)$ is equivalent
    to $(K_4,\emptyset)$ or the signed complete graph $K_5$ whose negative edges
    form a $5$-cycle,
  \item $\fln(K_n,\sig) = 3$ otherwise.
  \end{enumerate}
\end{theorem}

The result of~\cite{macajova_nowhere-zero_2015} concerning signed complete
bipartite graphs involves a special family of graphs. The family $\TT$
contains all signed bipartite graphs $(K_{3,n},\sig)$ such that
$n\geq 3$ and each vertex is incident with at most one negative edge,
except possibly for one vertex $v$ of degree $n$, which is incident
with at least one and at most $n-2$ negative edges; moreover, the two
vertices of degree $n$ other than $v$ are incident with $0$ and $1$
negative edges, respectively.

\begin{theorem}\label{t:complete-bi}
  Let $(K_{m,n},\sig)$ be a flow-admissible signed complete bipartite
  graph. Then
  \begin{enumerate}[label=(\roman*)]
  \item $\fln(K_{m,n},\sig) = 2$ if and only if $m,n,\size\sig$ are even,
  \item $\fln(K_{m,n},\sig) = 4$ if and only if either
    $(K_{m,n},\sig)\in\TT$, or $m=n=4$ and $\size\sig$ is odd,
  \item $\fln(K_{m,n},\sig) = 3$ otherwise.
  \end{enumerate}
\end{theorem}
\smallskip

Kaiser and Rollov\'{a}~\cite{kaiser_nowhere-zero_2014} studied
nowhere-zero flows in signed series-parallel graphs and showed that
Bouchet's Conjecture holds in this class:
\begin{theorem}
  Every flow-admissible signed series-parallel graph admits a
  nowhere-zero $6$-flow.
\end{theorem}
They pointed out that
although every unsigned series-parallel graph admits a nowhere-zero
$3$-flow and so these graphs are uninteresting from the point of view
of the $5$-Flow Conjecture, there are signed series-parallel graphs
with flow number $6$ (such as the one in Figure~\ref{fig:bouchet}(c)
and the other graphs in this family found by Schubert and Steffen in ~\cite{schubert_nowhere-zero_2015}).
\smallskip

Another class of graphs for which Bouchet's Conjecture has been
verified (see~\cite{schubert_nowhere-zero_2015}) is that of \emph{Kotzig
  graphs}, namely $3$-edge-colourable cubic graphs $G$ with the
property that the removal of the edges of any one colour produces a
Hamiltonian circuit of $G$.
\smallskip

A different approach to Bouchet's Conjecture was given by Rollov\' a et al. in~\cite{rollova_signed_2016}. Consider a flow-admissible signed graph $(G,\sig)$ with a minimum signature (that is a representative of the class of signed graphs equivalent to $(G,\sig)$ with minimum number of negative edges). If $|\sig|=0$, then the signed graph is all-positive and by Seymour's 6-flow theorem~\cite{seymour_nowhere-zero_1981}, it admits a nowhere-zero 6-flow. By Corollary~\ref{t:one-edge}, $|\sig|\neq 1$. Thus, the least
number of negative edges for which Bouchet's Conjecture is open is $2$. In~\cite{rollova_signed_2016} it was proved that the flow number of $(G,\sig)$ with $|\sig|=2$ is at most $7$, and if Tutte's 5-flow conjecture holds, then it is $6$. Moreover, Bouchet's Conjecture is true for cubic $(G,\sig)$ with $|\sig|=2$ if $G$ contains a bridge or if $G$ is $3$-edge-colourable or if $G$ is a critical snark (that is a graph that does not admits a nowhere-zero $4$-flow such that $G-e$ admits a nowhere-zero $4$-flow for any edge $e$ of $G$).  Furthermore, if $G$ is bipartite, then $(G,\sig)$ admits a nowhere-zero $4$-flow.  

\subsection{Circular flows}
\label{sec:circular}

Raspaud and Zhu~\cite{raspaud_circular_2011} generalised the concept
of unsigned circular flow to signed graphs. For $t\geq 2$ real, a
\emph{circular $t$-flow} in $(G,\sig)$ is an $\RR$-flow $(D,\varphi)$
such that for each edge $e$,
\begin{equation*}
  1\leq |\varphi(e)| \leq t-1.
\end{equation*}
If there is $t$ such that $(G,\sig)$ has a circular $t$-flow, then its
\emph{circular flow number} $\flnc(G,\sig)$ is the infimum of all such
$t$. The circular flow number is infinite if $(G,\sig)$ admits no
circular $t$-flow for any $t$.

\begin{figure}
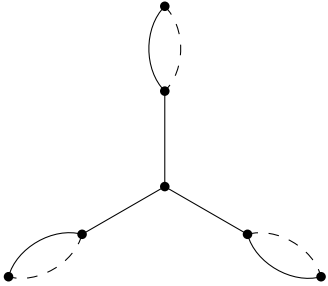

  \centering
  \fig{10}
  \caption{One member of a family of signed graphs for which the
    difference $\fln-\flnc$ tends to $2$.  }
  \label{fig:schubert}
\end{figure}

It is natural to ask whether the relation between $\fln$ and $\flnc$
for signed graphs is as close as in the unsigned case, where it is
known that $\fln(G) = \lceil\flnc(G)\rceil$ for any unsigned graph
$G$ --- in other words, the difference between $\fln$ and $\flnc$ is
less than $1$. Raspaud and Zhu~\cite{raspaud_circular_2011}
established the following upper bound on $\fln$:
\begin{theorem}[Raspaud and Zhu]\label{t:flnc-upper}
  For any flow-admissible signed graph $(G,\sig)$,
  \begin{equation*}
    \fln(G,\sig) \leq 2\lceil\flnc(G,\sig)\rceil-1.
  \end{equation*}
\end{theorem}
Furthermore, it was conjectured in~\cite{raspaud_circular_2011} that
$\fln(G,\sig)-\flnc(G,\sig)<1$ just as in the unsigned case.
However, the conjecture was disproved by Schubert and
Steffen~\cite{schubert_nowhere-zero_2015} who constructed signed
graphs $(G,\sig)$ with the difference $\fln(G,\sig)-\flnc(G,\sig)$
arbitrarily close to $2$. (See Figure~\ref{fig:schubert} for one
element of a family mentioned in~\cite{schubert_nowhere-zero_2015}.)
This result was strenthened by
\Macajova{} and Steffen~\cite{macajova_difference_2015} who proved that the
difference $\fln(G,\sig)-\flnc(G,\sig)$ can be arbitrarily close
to $3$. Furthermore, the result of~\cite{macajova_difference_2015}
shows that Theorem~\ref{t:flnc-upper} cannot be improved in general.

Bounds for the circular flow number of sufficiently edge-connected
signed graphs have been obtained by Raspaud and
Zhu~\cite{raspaud_circular_2011}:
\begin{theorem}[Raspaud and Zhu]\label{t:rz}
  Let $(G,\sig)$ be a flow-admissible signed graph. If $G$ is $6$-edge-connected, then $\flnc(G,\sig) < 4$.
\end{theorem}

Thomassen~\cite{thomassen_weak_2012} showed that with increasing edge-connectivity of a (non-signed) graph, the circular flow number approaches $2$.  Zhu~\cite{zhu_circular_2015} studied the existence of circular $t$-flows in signed graphs for $t$ close to $2$. 
It turns out that that it is not the case with signed graphs.
It is easy to see~\cite{zhu_circular_2015} that a graph that has exactly $2k+1$ negative edges has circuilar flow number at least $2+1/k$. Therefore it is natural to introduce
the notion of \emph{essentially $(2k+1)$-imbalanced} signed graph, defined as a signed graph such that every equivalent signed graph has either at least $2k+1$ negative
edges, or an even number of negative edges. 
The analogue to the result of Thomassen~\cite{thomassen_weak_2012}, shows that if both edge-connecticvity and essential imbalancednes rises, then the circular flow number goes to $2$~\cite{zhu_circular_2015}.
\begin{theorem}[Zhu]\label{t:imbal}
  Let $k$ be a positive integer. If a signed graph is essentially
  $(2k+1)$-imbalanced and $(12k-1)$-edge-connected, then its circular
  flow number is at most $2+\frac1k$.
\end{theorem}

Another interesting aspect of circular flows is,
that there exist gaps among the values of circular flow number of
regular graphs \cite{steffen_circular_2001} that separate bipartite
graphs from non-bipartite graphs. The gaps remain even in signed
case~\cite{schubert_nowhere-zero_2015}:

\begin{theorem}
Let $t\geq 1$ be an integer and $(G,\sig)$ be a signed $(2t+1)$-regular graph. If  $\flnc(G,\sig)=r$, then $r=2+\frac{1}{t}$ or $r\geq 2+\frac{2}{2t-1}$.
\end{theorem}

\subsection{Zero-sum flows}
\label{sec:zero-sum}

Let $G$ be a graph. A \emph{zero-sum flow} in $G$ (also called
`unoriented flow') is a function $f:\,E(G)\to\RR-\Setx0$ such that for
each vertex $v$, the values $f(e)$ of all edges $e$ incident with $v$
sum to zero. For a positive integer $k$, a zero-sum flow $f$ is a
\emph{zero-sum $k$-flow} if for each edge $e$, $f(e)$ is an integer
and $1\leq|f(e)|\leq k-1$. The \emph{zero-sum flow number} of a graph $G$, which will be denoted
by $\flns(G)$, is defined similarly to the flow number, namely as the
least $k$ such that $G$ admits a zero-sum $k$-flow (possibly
$\infty$). 

Zero-sum flows are tightly connected with nowhere-zero flows on signed
graphs. Indeed, a zero-sum $k$-flow in $G$ is nothing but a
nowhere-zero $k$-flow in the all-negative graph
$(G,E(G))$. Conversely, we can find a nowhere-zero flow on $(G,\sig)$
by finding a zero-sum flow of graph $G'$ created from $(G,\sig)$ by
subdividing all positive edges with one vertex of degree
$2$~\cite{akbari_zero-sum_2009}.  The $6$-flow Conjecture, therefore,
naturally translates to the following
conjecture~\cite{akbari_zero-sum_2009}:
\begin{conjecture}[Zero-Sum Conjecture]\label{conj:zsc}
  If a graph $G$ admits a zero-sum flow, then it admits a zero-sum
  $6$-flow.
\end{conjecture}





Much of the work regarding zero-sum flows pertains to regular
graphs. While not all $2$-regular graphs admit a zero-sum flow, Akbari
et al.~\cite{akbari_zero-sum_2009} proved that any $r$-regular graph
with $r\geq 4$ even has a zero-sum $3$-flow, and any cubic graph has a
zero-sum $5$-flow. This led to the following conjecture, stated
in~\cite{akbari_zero-sum_2010}:
\begin{conjecture}\label{conj:zs-regular}
  Any $r$-regular graph with $r\geq 3$ has a zero-sum $5$-flow.
\end{conjecture}

It was shown in~\cite{akbari_zero-sum_2010} that
Conjecture~\ref{conj:zs-regular} holds for odd $r$ divisible by $3$,
and that $\flns(G)\leq 7$ for general $r$. This was improved
in~\cite{akbari_note_2012} where the conjecture was proved for all $r$
except $r=5$. Furthermore, a mention in~\cite{sarkis_zero-sum_2015}
indicates that the case of $5$-regular graphs has also been settled in
the affirmative, namely in~\cite{zare_nowhere-zero_2013}. If so, then
Conjecture~\ref{conj:zs-regular} is proved.
\smallskip

Wang and Hu~\cite{wang_zero-sum_2012} gave exact values of $\flns(G)$
for certain $r$-regular graphs $G$:
\begin{itemize}
\item $\flns(G) = 2$ if either $r\equiv 0\pmod 4$, or
  $r\equiv 2\pmod 4$ and $\size{E(G)}$ is
  even,
\item $\flns(G) = 3$ if one of the following holds:
  \begin{itemize}[label=--]
  \item $r\equiv 2\pmod 4$ and $\size{E(G)}$ is odd,
  \item $r\geq 3$ is odd and $G$ has a perfect matching,
  \item $r\geq 7$ is odd and $G$ is bridgeless.
  \end{itemize}
\end{itemize}
The zero-sum flow number is determined in~\cite{wang_zero-sum_2012}
for the Cartesian product of an $r$-regular graph $G$ with an
$s$-vertex path $P_s$ ($r,s\geq 2$); namely, $\flns(G\times P_s) = 2$
if $r$ is odd and $s=2$, and $\flns(G\times P_s) = 3$ in the other
cases.

Further results concerning the zero-sum flow number of (not
necessarily regular) graphs include:
\begin{itemize}
\item $\flns(G)\leq 6$ if $G$ is bipartite
  bridgeless~\cite{akbari_zero-sum_2009},
\item $\flns(G)\leq 12$ if $G$ is
  hamiltonian~\cite{akbari_nowhere-zero_2015}.
\end{itemize}

As for particular graph classes, exact values of the zero-sum flow
number are known for graphs of the following types:
\begin{itemize}
\item wheels and fans~\cite{wang_zero-sum_2014},
\item all induced subgraphs of the hexagonal grid~\cite{wang_zero-sum_2013},
\item certain subgraphs of the triangular
  grid~\cite{wang_zero-sum_2014}.
\end{itemize}
 
Results on the computational complexity of deciding the existence of a
zero-sum flow were obtained in~\cite{dehghan_complexity_2015}. For
positive integers $k\leq\ell$, let us define a \emph{$(k,\ell)$-graph}
to be an unsigned graph with minimum degree at least $k$ and maximum
degree at most $\ell$. It is known~\cite{akbari_zero-sum_2009} that
Conjecture~\ref{conj:zsc} is equivalent to its restriction to
$(2,3)$-graphs (that is, subdivisions of cubic graphs).

Let us summarise the findings of~\cite{dehghan_complexity_2015}:
\begin{theorem}
  \begin{enumerate}[label=(\roman*)]
  \item There is a polynomial-time algorithm to determine whether a
    given $(2,4)$-graph with $O(\log n)$ vertices of degree $4$ has
    a zero-sum $3$-flow,
  \item it is NP-complete to decide whether a given $(3,4)$-graph
    admits a zero-sum $3$-flow,
  \item it is NP-complete to decide whether a given $(2,3)$-graph
    admits a zero-sum $4$-flow.
  \end{enumerate}
\end{theorem}

In~\cite{sarkis_zero-sum_2015}, zero-sum flows are viewed in the more
general setting of matrices. A \emph{zero-sum $k$-flow} of a real
matrix $M$ is an integer-valued vector in its nullspace with nonzero
entries of absolute value less than $k$. From this point of view, a
zero-sum flow in a graph is a zero-sum flow of its vertex-edge
incidence matrix. The main result of~\cite{sarkis_zero-sum_2015}
proves the existence of a zero-sum $k$-flow of the incidence matrix
for linear subspaces of dimension $1$ and $m$, respectively, of the
vector space $\FF_q^n$ (where $q$ is a prime power, $0\leq m\leq n$
and $k$ is a suitable integer depending on $q,m$ and $n$).


\subsection{Structural results}
\label{sec:structural}

Zaslavsky~\cite{zaslavsky_signed_1982-1} proved that any signed graph
$(G,\sig)$ is associated with a matroid, called the \emph{frame
  matroid} of $(G,\sig)$ and denoted by $M(G,\sig)$. (For background
on matroid theory, see, e.g.,~\cite{oxley_matroid_1992}.) One way to
define this matroid is to specify its circuits, and these are
precisely the signed circuits introduced in
Section~\ref{sec:intro}. Properties of $M(G,\sig)$ are studied
in~\cite{zaslavsky_signed_1982-1} and subsequently in a number of
papers. See, for
instance,~\cite{slilaty_decompositions_2007,slilaty_projective-planar_2007}
(decomposition theorems), \cite{gerards_graphs_1990,qin_regular_2009}
(minors) and the references therein.

Let us call a $\ZZ$-flow in a (signed or unsigned) graph
\emph{nonnegative} if all its values are nonnegative. A nonnegative
$\ZZ$-flow is \emph{irreducible} if it cannot be expressed as a sum of
two nonnegative $\ZZ$-flows. Tutte~\cite{tutte_class_1956} showed that
in an unsigned graph, irreducible nonnegative $\ZZ$-flows are
precisely the \emph{circuit flows} (elementary flows along
circuits). Equivalently, any nonnegative $\ZZ$-flow is a sum of
circuit flows. 

\begin{figure}
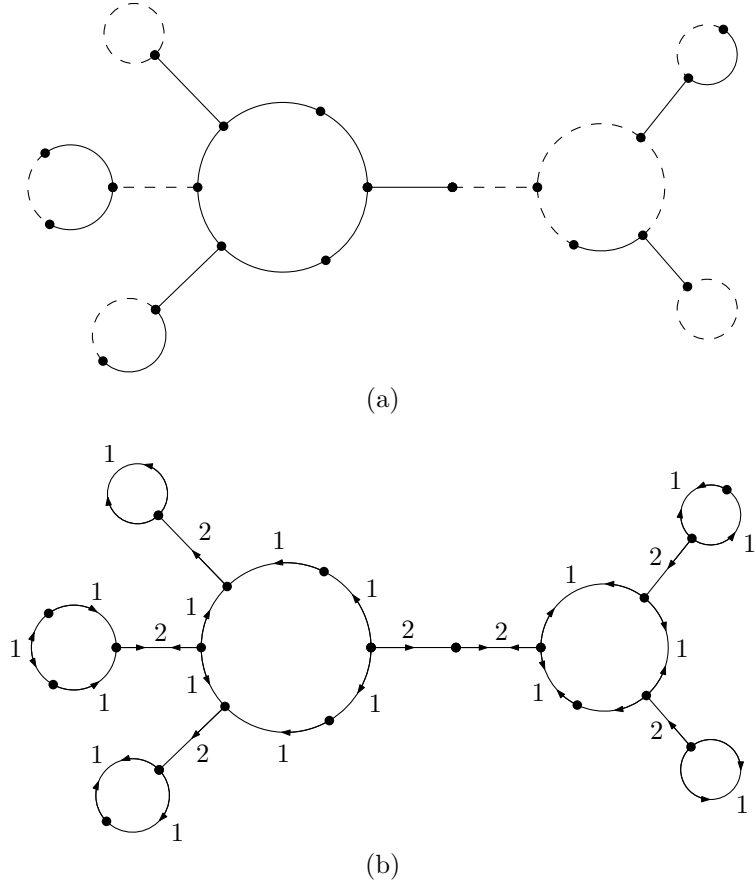

  \centering
  \subfloat[][]{\fig{11}}\\
  \subfloat[][]{\fig{12}}  
  \caption{(a) An essential signed graph $(H,\Theta)$. (b) An essential
    flow in $(H,\Theta)$.}
  \label{fig:essential}
\end{figure}

For signed graphs, a characterisation of irreducible flows appears in
the (apparently unpublished) manuscript~\cite{chen_resolution_2007}
and in~\cite{chen_classification_2011}. Let us say that a signed graph
$(H,\Theta)$ is \emph{essential} if it satisfies the following
properties:
\begin{itemize}
\item the degree of any vertex of $H$ is $2$ or $3$,
\item each vertex is contained in at most one circuit of $H$,
\item for any circuit $C$ of $H$, $\size{E(C)\cap\Theta}$ has the same
  parity as the number of bridges of $H$, incident with $C$.
\end{itemize}
(See Figure~\ref{fig:essential} for an example.)

Let $C$ be a circuit in $H$ and let $\Theta' = \Theta\cap E(C)$. An
orientation $D$ of $(C,\Theta')$ is \emph{consistent} at a vertex $v$
of $C$ if $v$ is incident with precisely one incoming
half-edge. Otherwise, $D$ is \emph{inconsistent} at $v$. A vertex is a
\emph{source} (\emph{sink}) in an orientation of a signed graph if all
incident half-edges are outgoing (incoming, respectively).

Given an essential signed graph $(H,\Theta)$, an \emph{essential flow}
in $(H,\Theta)$ is obtained as follows:
\begin{itemize}
\item choose an orientation of $(H,\Theta)$ with no sources nor sinks,
  such that the induced orientation of any circuit $C$ of $H$ is
  inconsistent at any vertex incident with a bridge of $H$ (we call
  this an \emph{essential orientation of $(H,\Theta)$}),
\item assign flow value $2$ to each bridge of $H$ and value $1$ to the
  other edges.
\end{itemize}

The main result of~\cite{chen_resolution_2007} reads as follows:
\begin{theorem}
  Irreducible flows in a signed graph $(G,\sig)$ are exactly the
  essential flows on essential subgraphs of $(G,\sig)$.
\end{theorem}

Independently, \Macajova{} and
\Skoviera~\cite{macajova_characteristic_2014} proved the closely
related result that nonnegative flows in signed graphs can be obtained
as sums of elementary flows along signed circuits at the cost of
allowing half-integral values. Let $C$ be a signed circuit in
$(G,\sig)$. An orientation of $C$ will be called \emph{proper} if it
is an essential orientation of the essential subgraph $C$ as defined
above. The \emph{characteristic flow} of $C$ is the function
$\chi_C:\,E(G)\to\QQ$ with values
\begin{equation*}
  \chi_C(e) =
  \begin{cases}
    1 & \text{if $C$ is balanced and $e\in E(C)$},\\
    \frac12 & \text{if $C$ is unbalanced and $e$ is contained in an ordinary
    circuit of $C$,}\\
    1 & \text{if $C$ is unbalanced and $e$ is contained in the
      connecting path of $C$,}\\
    0 & \text{otherwise.}
  \end{cases}
\end{equation*}
\begin{theorem}\label{t:sum}
  Let $(D,\varphi)$ be a $\ZZ$-flow on a signed graph $(G,\sig)$ such
  that all values of $\varphi$ are nonnegative. Then $D$ contains
  proper orientations of signed circuits $C_1,\dots,C_m$ and positive
  integers $\alpha_i$ ($i=1,\dots,m$) such that
  \begin{equation*}
    \varphi = \sum_{i=1}^m \alpha_i \chi_{C_i}.
  \end{equation*}
\end{theorem}

Chen and Wang~\cite{chen_flow_2009} defined cuts in signed graphs and
introduced the circuit and bond lattices of signed graphs.


\section{Circuit covers}
\label{sec:circuit-covers}

The Shortest Cycle Cover Conjecture of Alon and
Tarsi~\cite{alon_covering_1985} asserts that the edge set of any
(unsigned) bridgeless graph with $m$ edges can be covered with
circuits of total length at most $\tfrac{7m}5$. For background on the
conjecture and its relation to other problems see,
e.g.,~\cite[Chapter~14]{zhang_circuit_2012}. The best general upper
bound to the Shortest Cycle Cover Conjecture is due to Bermond et
al.~\cite{bermond_shortest_1983}, and Alon and
Tarsi~\cite{alon_covering_1985}, who independently proved that any
bridgeless graph admits a circuit cover of total length at most
$\tfrac{5m}3$. 

The signed version of the problem was introduced by \Macajova{} et
al.~\cite{macajova_circuit_2016}. In a signed graph $(G,\sig)$, it is
natural to seek a cover of the entire edge set by signed circuits (a
\emph{signed circuit cover}) of small length. Here, the length of a
signed circuit is the number of its edges, and the \emph{length} of
the cover is defined as the sum of lengths of its elements. Similarly as in the case of flows,
for a signed circuit cover of $G$ to exist, each edge has to be contained in a signed circuit, which is equivalent
to $(G,\sig)$ being flow-admissible (see Theorem~\ref{t:flow-admissible}). The authors
of~\cite{macajova_circuit_2016} obtained the following bounds:
\begin{theorem}
  Let $(G,\sig)$ be a flow-admissible signed graph with $m$
  edges. Then there is a signed circuit cover of $(G,\sig)$ of length
  at most $11m$. Moreover, if $G$ is bridgeless, then it admits a
  signed circuit cover of length $9m$.
\end{theorem}

An application of Theorem~\ref{t:sum} given
in~\cite{macajova_characteristic_2014} yields the following result:
\begin{theorem}
  Let $(G,\sig)$ be a signed graph with $m$ edges. If $(G,\sig)$
  admits a nowhere-zero $k$-flow, then it has a signed circuit cover
  of length at most $2(k-1)m$.
\end{theorem}

Currently the best bound for this problem is due to Cheng et
al.~\cite{cheng_shortest_2015}. Let the \emph{negativeness} of a
signed graph be the minimum number of negative edges in any
equivalent signed graph.
\begin{theorem}\label{t:scc}
  Let $(G,\sig)$ be a flow-admissible signed graph with $n$ vertices,
  $m$ edges and negativeness $\varepsilon>0$. Then the following holds:
  \begin{enumerate}[label=(\roman*)]
  \item $(G,\sig)$ has a signed circuit cover of length at most
    \begin{equation*}
      m+3n+\min\Setx{\frac{2m}3+\frac{4\varepsilon}3-7,n+2\varepsilon-8}
      \leq \frac{14m}3-\frac{5\varepsilon}3-4.
    \end{equation*}
  \item If $G$ is bridgeless and $\varepsilon$ is even, then
  $(G,\sig)$ admits a signed circuit cover of length at most
  \begin{equation*}
    m+2n+\min\Setx{\frac{2m}3+\frac{\varepsilon}3-4,n+\varepsilon-5}.
  \end{equation*}
  \end{enumerate}
\end{theorem}


\section{Relation to colouring}
\label{sec:relation-colouring}

As indicated in Section~\ref{sec:intro}, nowhere-zero flows on signed
graphs are dual to local tensions on non-orientable surfaces
(see~\cite{devos_flows_2014}). For an embedded directed graph $G$ and
a function $\phi: E(G) \to A$, we say that $\phi$ is a \emph{local
  tension} if the height of every facial walk is $0$ (the \emph{height} of a facial walk for a given orientation of the surface, is the sum of values on edges of the face directed consistently with the orientation minus the sum of values on the edges directed inconsistently). Since the
surface is non-orientable, the dual graph $G^*$ of a directed graph
$G$ will be bidirected, and the condition on the facial walk of a
local tension of $G$ will translate to the conservation law at a
vertex of $G^*$. This duality can be used to show that the signed
Petersen graph of Figure~\ref{fig:bouchet}(a) does not admit any
nowhere-zero $5$-flow, because its dual in the projective plane is
$K_6$, which does not admit any $5$-local-tension.
\smallskip

The flow-colouring duality does not hold for signed graphs for any of
the known definitions of signed colourings. However, given the
importance of the duality in the unsigned case, it will perhaps be
useful to include an overview of the essential definitions of signed
colourings and the appropriate literature for completeness. The first
definition is due to Zaslavsky~\cite{zaslavsky_signed_1982}:
\smallskip

\textbf{Signed colouring according to Zaslavsky.} If $(G,\sig)$ is a
signed graph, a (signed) colouring of $(G,\sig)$ in $k$ colors, or in
$2k+1$ signed colors is a mapping
$$f:\ V(G) \to \{-k,-k+1,\ldots, -1,0,1,\ldots, k-1,k\}.$$ A colouring
is zero-free if it never uses the value 0. The (signed) colouring is
proper if for every edge uv $f(u)\neq f(v)\cdot \sig(uv)$. Zaslavsky
further defined the chromatic polynomial $\xi_G(2k+1)$ of a signed
graph G to be the function counting the number of proper signed
colourings of $G$ in $k$ colours. The balanced chromatic polynomial
$\xi^b_G(2k)$ is the function which counts the zero-free proper signed
colourings in $k$ colours. Then the chromatic number $\gamma(G)$ of
$G$, according to Zaslavsky, is the smallest nonnegative integer $k$
for which $\xi_G(2k+1) > 0$. Furthermore, the strict chromatic number
$\gamma^*(G)$ of $G$ is the smallest nonnegative integer $k$ such that
$\xi^b_G(2k) > 0$. Neither of Zaslavsky’s definitions of chromatic
number is a direct extension of the usual chromatic number of an
unsigned graph. These signed colourings were further investigated by
Zaslavsky in the series of papers~\cite{zaslavsky_chromatic_1982,
  zaslavsky_how_1984, zaslavsky_signed_1995}, and recently by Beck et
al.~\cite{beck_chromatic_2015}, by Davis~\cite{davis_unlabeled_2015},
and by Schweser et al.~\cite{schweser_degree_2015}.
\smallskip

\textbf{Signed colouring according to \Macajova{}, Raspaud and \v Skoviera.} A
modification of Zaslavsky's definition of a signed colouring rised to
a 'similar' definition of signed colouring introduced by \Macajova{},
Raspaud,\v Skoviera in~\cite{macajova_chromatic_2014}. Let, for each
$n \geq 1$, $M_n\subseteq \ZZ$ be a set
$M_n=\{\pm 1, \pm 2, \ldots, \pm k\}$ if $n = 2k$, and
$M_n=\{0,\pm 1, \pm 2, \ldots, \pm k\}$ if $n = 2k + 1$. A proper
colouring of a signed graph $G$ that uses colours from $M_n$ is an
$n$-colouring. Thus, an $n$-colouring of a signed graph uses at most
$n$ distinct colours. Note that if $G$ admits an $n$-colouring, then
it also admits an $m$-colouring for each $m \geq n$. The smallest $n$
such that $G$ admits an $n$-colouring is the signed chromatic number
of $G$ and is denoted by $\xi(G)$. It is easy to see that the
chromatic number of a balanced signed graph coincides with the
chromatic number of its underlying unsigned graph, and hence this
definition of chromatic number differs from Zaslavsky's. Moreover,
$\xi(G) = \gamma(G) +\gamma^*(G)$. Despite the fact that this
definition is very recent, it has already been studied by Jin et
al.~\cite{jin_choosability_2016} and by Fleiner et
al.~\cite{fleiner_coloring_2016}.
\smallskip

\textbf{Signed colourings via signed homomorphisms according to
  Naserasr, Rollov\' a and Sopena.} A chromatic number for a graph can
be defined using graph homomorphism. This is the case also for signed
graphs when we use signed homomorphism.  The notion of homomorphism of
signed graphs was introduced by Guenin~\cite{guenin_packing_2005}, and
later systematically studied by Naserasr et
al.~\cite{naserasr_homomorphisms_2015} (cf. also the extended
abstract~\cite{naserasr_homomorphisms_2013-1}). Given two signed
graphs $(G,\sig_1)$ and $(H,\sig_2)$, we say there is a homomorphism
of $(G,\sig_1)$ to $(H,\sig_2)$ if there is $(G,\sig'_1)$, which is
equivalent to $(G,\sig_1)$, $(H,\sig'_2)$, which is equivalent to
$(H,\sig'_2)$, and a mapping $\Phi:\ V(G) \to V(H)$ such that every
edge of $(G,\sig'_1)$ is mapped to an edge of $(H,\sig'_2)$ of the
same sign. A signed chromatic number of a signed graph $(G,\sig)$ is
the smallest order of a signed graph to which $(G,\sig)$ admits a
homomorphism. The signed chromatic number of a balanced signed graph
is the same as the chromatic number of an unsigned graph, hence this
definition differs from the one by Zaslavsky. Moreover, it also
differs from the one by M\' a\v cajov\' a, Raspaud and \v Skoviera,
for example, for an unbalanced $C_4$ (the chromatic number of
unbalanced $C_4$ is 4 in terms of homomorphism, while it is 3 in the
other case). Signed homomorphisms were investigated by Naserasr et
al.~\cite{naserasr_homomorphisms_2013} (cf. also the extended
abstract~\cite{naserasr_homomorphisms_2013-2}), Foucaud et
al.~\cite{foucaud_complexity_2014}, Brewster et
al.~\cite{brewster_complexity_2015}, Ochem et
al.~\cite{ochem_homomorphisms_2014}, and recently by Das et
al. in~\cite{das_chromatic_2016} and Das et
al. in~\cite{das_relative_2016}.

\printbibliography


\end{document}
